\documentclass[10pt]{article}

\usepackage{amsmath}
\usepackage{amsfonts}

\usepackage{cite} 
\usepackage{url}  
\usepackage{ifthen}  
\usepackage{multicol}   
\usepackage[utf8]{inputenc} 
\def\includegraphics{}
\newboolean{publ}

\newenvironment{bmcformat}{\begin{raggedright}\baselineskip20pt\sloppy\setboolean{publ}{false}}{\end{raggedright}\baselineskip20pt\sloppy}

\newtheorem{Theorem}{Theorem}

\DeclareMathOperator{\e}{e}

\begin{document}
\begin{bmcformat}

\title{A new semi-analytical approach for numerical solving of Cauchy problem for functional differential equations}

\author{Josef Rebenda 
         \thanks{CEITEC BUT, Brno University of Technology, Technicka 3058/10, 61600 Brno, Czech Republic (josef.rebenda@ceitec.vutbr.cz)}
      \and
         Zden\v{e}k \v{S}marda 
         \thanks{Department of Mathematics, Faculty of Electrical Engineering and Communication, Brno University of Technology, Technicka 8, 616 00 Brno, Czech Republic (smarda@feec.vutbr.cz) - corresponding author}
       \and 
       Yasir Khan 
       \thanks{Department of Mathematics, Zhejiang University, Hangzhou 310027, China (yasirmath@yahoo.com)}%
       }



\date{}

\maketitle

\begin{abstract}
One of the major challenges of contemporary mathematics is numerical solving of various problems for functional differential equations (FDE), in particular Cauchy problem for delayed and neutral differential equations.
Recently large variety of methods to handle this task appeared.
In the paper, we present new semi-analytical approach for FDE's consisting in combination of the method of steps and a technique called differential transformation method (DTM).
This approach reduces the original Cauchy problem for delayed or neutral differential equation to Cauchy problem for ordinary differential equation for which DTM is convenient and efficient method.
Moreover, there is no need of any symbolic calculations or initial approximation guesstimates in contrast to methods like homotopy analysis method, homotopy perturbation method, variational iteration method or Adomian decomposition method.
The efficiency of the proposed method is shown on certain classes of FDE's with multiple constant delays including FDE of neutral type. We also compare it to the current approach of using DTM and the Adomian decomposition method where Cauchy problem is not well posed.
\end{abstract}

\ifthenelse{\boolean{publ}}{\begin{multicols}{2}}{}

\section*{Introduction}

 For the purpose of clarity, we consider the following functional differential equation of $n$-th order with multiple constant delays
\begin{equation}\label{1}
u^{(n)}(t) = f(t,u(t),u'(t),\dots,u^{(n-1)}(t),\mathbf{u}_1(t-\tau_1 ),\mathbf{u}_2(t-\tau_2 ), \dots, \mathbf{u}_r(t-\tau_r ) ),
\end{equation}
where $\mathbf{u}_i(t-\tau_i )= (u(t-\tau_i),u'(t-\tau_i),\dots,u^{(m_i)}(t-\tau_i))$ is $m_i$-dimensional vector function, $m_i \leq n$, $i=1,2,\dots,r$, $r \in \mathbb N$ and $f \colon [t_0,\infty) \times R^n \times R^\omega$ is a continuous function, where $\omega = \sum\limits_{i=1}^r m_i$.\\
Let $t^*= \max\{\tau_1,\tau_2,\dots,\tau_r\}$, $m= \max\{m_1,m_2,\dots,m_r\}$, $m \leq n $.  In case $m=n$ equation \eqref{1} is of neutral type, otherwise it is delayed differential equation. Initial function $\phi(t)$ needs to be assigned  for equation \eqref{1} on the interval $[t_0-t^*, t_0]$.
Furthermore, for the sake of simplicity, we assume that $\phi(t) \in C^n([t_0-t^*,t_0])$.

Investigation of equation \eqref{1} is important since there is plenty of applications of such equations in real life. As examples, we mention models describing behaviour of the central nervous system in a learning process, species populations struggling for a common food, dynamics of an autogenerator with delay and second-order filter, systems controlled by PI or PID regulators, evolution of population of one species etc. For further models and details, see e.g. \cite{kolmanovski}.

In the last two decades, various methods such as homotopy analysis method (HAM) \cite{wang}, \cite{alomari}, homotopy  perturbation method (HPM) \cite{shakeri}, variational iteration method (VIM) \cite{chen},  Adomian decomposition method (ADM)  \cite{evans}, \cite{cocom}, \cite{saeed},
Taylor polynomial method \cite{sezer}, Taylor collocation method \cite{bellour}
and differential transformation method (DTM) \cite{smarda}, \cite{arikoglu}, \cite{karakoc}, \cite{mohammed} have been considered to approximate solutions of certain classes of equation \eqref{1} in a series form. However, in several papers, for example Evans and Raslan \cite{evans}, Mohammed and Fadhel \cite{mohammed} initial problems are not properly defined. The authors use only initial conditions in certain points, not the initial function on the whole interval, thus the way to obtain solutions of illustrative examples is not correct. Moreover, transformation formulas used in the calculations are very complicated.

The crucial idea of our concept is to combine differential transformation method and general method of steps (more details on method of steps can be found for instance in monographies Kolmanovskii and Myshkis \cite{kolmanovski} or Bellen and Zennaro \cite{bellen}).
This approach enables us to replace the terms involving delay with initial function and its derivatives. Consequently, the original Cauchy problem for delayed or neutral differential equation is reduced to Cauchy problem for ordinary differential equation. Also the ambiguities mentioned above are removed. Further, while ADM, HAM, HPM and VIM require initial approximation guess and symbolic computation of necessary derivatives and, in general, $n$-dimensional integrals in iterative schemes, presented method is different: Cauchy problem for FDE is reduced to a system of recurrence algebraic relations.

\section{Preliminaries}
This section contains several well-known results which we compare to our approach later.

The concept of differential transformation
was used by Zhou \cite{zhou} who applied it to solve linear and nonlinear initial value problems in electrical circuit analysis.

Differential transformation of the $k-$th derivative of function $u(t)$ is defined as
\begin{equation}\label{2}
U(k) = \frac{1}{k!} \left[ \frac{d^ku(t)}{dt^k} \right]_{t=t_0},
\end{equation}
where $u(t)$ is the original function and $U(k)$ is the transformed function. Inverse differential transformation of $U(k)$ is defined as follows:
\begin{equation}\label{3}
u(t) = \sum_{k=0}^{\infty}U(k)(t-t_0)^k. 
\end{equation}

Many transformation formulas were derived from definitions \eqref{2} and \eqref{3} in the past. We mention several formulas for illustration:
\begin{Theorem}[\cite{smarda}, v) and vi) from definition] \label{t1}
Assume that $F(k)$, $G(k)$, $H(k)$ and $U_i(k)$, $i=1,\dots,n$, are differential transformations of functions $f(t)$, $g(t)$, $h(t)$ and
$u_i(t)$, $i=1,\dots,n$, respectively. Then
\begin{gather*}
\begin{array}{lllcl}
i) & {\rm  If} & f(t) = {\displaystyle \frac{d^ng(t)}{dt^n}}, & {\rm then} & F(k) = {\displaystyle \frac{(k+n)!}{k!}}G(k+n). \\[2mm]
ii)& {\rm  If} & f(t) = g(t)h(t), & {\rm then} & F(k) = \sum_{l=0}^k G(l)H(k-l). \\[2mm]
iii)& {\rm  If} & f(t) = t^n, & {\rm then} & F(k) =\delta (k-n),\ \delta\ {\rm is\ the\ Kronecker\ delta\  symbol}, t_0=0. \\[1mm]
iv)& {\rm  If} & f(t) = e^{\lambda t}, & {\rm then} & F(k) = {\displaystyle \frac{\lambda^k}{k!}}, t_0=0. \\[2mm]
v)& {\rm  If} & f(t) = \sin t,  & {\rm then} & F(k) =
 \begin{cases}
(-1)^{\frac{k-1}{2}} \frac{1}{k!} & \text{if } k=2n+1, n \in \mathbb{N}_0 \\
  0 & \text{if } k=2n, n \in \mathbb{N}_0
 \end{cases}, t_0=0. \\[2mm]
vi)& {\rm  If} & f(t) = \cos t,  & {\rm then} & F(k) =
\begin{cases}
(-1)^{\frac{k}{2}} \frac{1}{k!} & \text{if } k=2n, n \in \mathbb{N}_0 \\
  0 & \text{if } k=2n+1, n \in \mathbb{N}_0
 \end{cases}, t_0=0. \\[2mm]
vii) &{\rm  If} & f(t) = \prod_{i=1}^n u_i(t),& {\rm then}&
\end{array}
\\
F(k)= \sum_{s_1=0}^k \sum_{s_2=0}^{k-s_1} \dots \sum_{s_n=0}^{k-s_1-\cdots s_{n-1}} U_1(s_1)\dots U_{n-1}(s_{n-1}) U_n(k-s_1-\dots -s_n).
\end{gather*}
\end{Theorem} 
\begin{Theorem}[\cite{arikoglu}, \cite{karakoc}, \cite{mohammed}] \label{t2}
Assume that $F(k)$, $G(k)$ are differential transformations of functions $f(t)$, $g(t)$, where $a>0$ is a real constant. If
\begin{gather}
\begin{array}{lcl}
 f(t) = {\displaystyle g(t-a)}, & {\rm then} & {\displaystyle F(k) =  \sum_{i=k}^{N}(-1)^{i-k} \binom{i}{k} a^{i-k} G(i), \ N \rightarrow \infty.}
\end{array}
\end{gather} \label{4}
\end{Theorem}
Using Theorem \ref{t2} and the formula i) in Theorem \ref{t1}, differential transformation formula for function $f(t) = \frac{d^n}{dt^n}g(t-a)$ can be deduced.
\begin{Theorem}[\cite{arikoglu}] \label{t3}
Assume that $F(k)$, $G(k)$ are differential transformations of functions $f(t)$, $g(t)$, $a>0$. If
\begin{gather}
\begin{array}{lcl}
f(t) = {\displaystyle \frac{d^n}{dt^n}g(t-a)}, & {\rm then} & {\displaystyle F(k)= \frac{(k+n)!}{k!}\sum_{i=k+n}^{N}(-1)^{i-k-n} \binom{i}{k+n} a^{i-k-n}G(i), \ N \rightarrow \infty.}
\end{array}
\end{gather} \label{5}
\end{Theorem}
Using Theorems \ref{t2}, \ref{t3} and formula vii) in Theorem \ref{t1} any differential transformation of a product of functions with delayed arguments and derivatives of that functions can be proved. However,
such formulas are complicated and not easy applicable for solving functional differential equations with multiple constant delays (see for example \cite{arikoglu}, \cite{karakoc}, \cite{mohammed}).

\section{Main results}
Consider equation \eqref{1} subject to initials conditions
\begin{equation}\label{6}
u(t_0)=u_0, u'(t_0)=u_1, \dots, u^{(n-1)}(t_0) = u_{n-1}
\end{equation}
and subject to initial function $\phi(t)$ on interval $[t_0-t^*, t_0]$ such that
\begin{equation}\label{6'}
\phi(t_0) = u(t_0), \phi'(t_0)= u'(t_0), \dots, \phi^{(n-1)}(t_0)= u^{(n-1)}(t_0).
\end{equation}
First we apply the method of steps. We substitute the initial function $\phi(t)$ and its derivatives in all places where unknown functions with deviating arguments and derivatives of that functions appear. Then equation \eqref{1} changes to ordinary differential equation
\begin{equation}\label{7}
u^{(n)}(t) = f(t,u(t),u'(t),\dots,u^{(n-1)}(t),\mathbf{\Phi}_1(t-\tau_1),\mathbf{\Phi}_2(t-\tau_2), \dots, \mathbf{\Phi}_r
(t-\tau_r) ),
\end{equation}
where $\mathbf{\Phi}_i(t-\tau_i )= (\phi(t-\tau_i),\phi'(t-\tau_i),\dots,\phi^{(m_i)}(t-\tau_i))$, $m_i \leq n$, $i=1,2,\dots,r$.
Now applying DTM we get recurrence  equation
\begin{equation}\label{8}
\frac{(k+n)!}{k!} U(k+n) = \mathcal F \Bigl( k, U(0), U(1), \dots, U(k+n-1)\Bigr).
\end{equation}
Using transformed initial conditions and then inverse transformation rule, we obtain approximate solution of
equation \eqref{1} in the form of infinite Taylor series
$$
u(t) = \sum_{k=0}^{\infty} U(k)(t-t_0)^k
$$
on the interval $[t_0, t_0+ \alpha]$, where $\alpha = \min\{\tau_1,\tau_2, \dots, \tau_r\}$, and $ u(t) = \phi(t)$ on the interval $[t_0-t^*,t_0]$.
We demonstrate potentiality of our approach on several examples.

{\bf Example 1.} Consider the problem that was solved using Taylor series method by Sezer and Akyuz-Dascioglu \cite{sezer} and using DTM by Arikoglu and Ozkol \cite{arikoglu},
\begin{equation}\label{9}
u''(t) -tu'(t-1) +u(t-2) =-t^2 -2t + 5, 
\end{equation}
\begin{equation}\label{10}
 u(0)=-1, \ u'(-1) = -2, \quad -2 \leq t \leq 0.
 \end{equation}
First, we remark that such formulation of problem is not correct since if we take for instance $t=-1$, then $u(t-2) = u(-3)$ is not defined at all. If we omit condition $-2 \leq t \leq 0$ in \eqref{10}, then according to Sezer and Akyuz-Dascioglu \cite{sezer} we are looking for solutions of a problem involving mixed conditions. However, this is not a Cauchy problem, thus it is not clear what kind of solution are we looking for since uniqueness of solution is not guaranteed. Let us see consequence of this fact.\\
Applying  current approach of using DTM  on both  sides of equation \eqref{9} and conditions \eqref{10}, Arikoglu and Ozkol \cite{arikoglu} obtained  recurrence relation
\begin{eqnarray}\label{11}
&\ &(k+1)(k+2)U(k+2) - \sum_{k_1=0}^k \sum_{h_1=k-k_1+1}^N (k-k_1+1)\binom{h_1}{k-k_1+1} (-1)^{h_1-k+k_1-1} U(h_1)\delta(k_1-1) \nonumber \\
&+& \sum_{h_1=k}^N \binom{h_1}{k} (-2)^{h_1-k}U(h_1) = -\delta(k-2) -2\delta(k-1) + 5\delta(k) 
\end{eqnarray}
 and transformed mixed conditions
$$
U(0)=-1, \quad \sum_{k=0}^N kU(k) (-1)^{k-1} = -2 .
$$
Taking $N=4$ and using inverse differential transformation formula \eqref{3}, they claimed solution $u(t) = -1+t^2$  which has the same form for all $N>4$. It can be easily verified that this $u(t)$ is exact solution of \eqref{9}, \eqref{10} for all $t \in \mathbb R$. Nevertheless, this coincidence is possible only for solutions in expected polynomial form. Generally, since $N$ is a finite number, we only get  approximate solutions.\\
Applying different approach using Taylor series method, Sezer and Akyuz-Dascioglu \cite{sezer} obtained one-parameter class of solutions in the form $u(t)=t^2-1+a_1 (t+2)$, where $a_1 \in \mathbb R$ is a parameter. It is easy to check that, again, such $u(t)$ is exact solution of \eqref{9}, \eqref{10} for all $t \in \mathbb R$. Furthermore, this class of solutions contains the solution achieved by Arikoglu and Ozkol \cite{arikoglu}.

Now recall that to formulate Cauchy problem for \eqref{9} correctly, we have to prescribe initial function $\phi (t)$ on $[-2,0]$ satisfying conditions \eqref{6'} as well:
\begin{equation}\label{12}
u(t)=\phi (t) \ \text{for } t \in [-2,0), \ u(0)= \phi (0) = u_0, \ u'(0) = \phi' (0) = u_1.
\end{equation}
Thus a solution of Cauchy problem \eqref{9}, \eqref{12} then should be expected in the form
\begin{equation*}
u(t)=
\begin{cases}
   \psi(t), \ t \in [0,1] \\
   \phi (t), \  t \in [-2,0]
 \end{cases}.
\end{equation*}
It is obvious that neither current approach of using DTM nor the other approach is suitable for solving Cauchy problem \eqref{9}, \eqref{12} since, in general, there are infinitely many initial functions satisfying conditions \eqref{10} (see \eqref{13} below as an example), whereas those approaches always give the same result.

On the other hand, our approach enables to find unique solution of Cauchy problem \eqref{9}, \eqref{12}. For initial function $\phi(t)= -1+t^2$ satisfying \eqref{10} we have simple recurrence relation
$$
U(k+2) = \frac{2\delta(k)}{(k+1)(k+2)}, \ k\geq 0.
$$
From this relation and the fact that $U(0)=-1$, $U(1)=0$, we get $U(2)= 1$, $U(k) = 0$ for $ k \geq3$. Hence the exact solution is $u(t)= -1+t^2$, generally on $[0,\infty)$, and this solution of Cauchy problem \eqref{9}, \eqref{12} is unique.

Further, for initial function 
\begin{equation} \label{13}
\phi(t) = -t^2 -4t-1 
\end{equation}
that also satisfies conditions \eqref{10}, using DTM combined with the method of steps we obtain different solution of \eqref{9}, \eqref{12}. In this case we have recurrence relation
\begin{equation}\label{14}
U(k+2)= \frac{-2\delta(k-2) -4\delta(k-1) +2\delta(k)}{(k+1)(k+2)}, \ k\geq 0.
\end{equation}
From initial conditions $u(0)=0$, $u'(0)=-4$, we have $U(0)=0$, $U(1)=-4$ and from \eqref{14} we get $U(2)=1$, $U(3)=-2/3$, $U(4) = -1/6$, $U(k) = 0$ for $k \geq 5$. 
Thus we have exact solution in the form
\begin{equation*}
u(t)=
\begin{cases}
  -1-4t+t^2 -\frac{2}{3}t^3 -\frac{1}{6}t^4,  \ t \in [0,1]\\
                 -t^2-4t-1, \  t \in [-2,0]
 \end{cases},
\end{equation*}
and, again, this solution is unique. However, it is not possible to obtain this solution using either current DTM approach or the other mentioned method.\\
{\bf Example 2.}
Consider Cauchy problem consisting of first order differential equation with one constant delay
\begin{equation} \label{15}
u' (t) = \frac{1}{a} u(t) - \frac{1}{a} u(t-a) +a, \ a>0,
\end{equation}
and initial function
\begin{equation} \label{16}
\phi (t) =t^2, \ t \in [-a,0].
\end{equation}
Obviously, $u(t) = t^2$  is unique solution of Cauchy problem \eqref{15}, \eqref{16} on $[0,\infty)$.
From initial function \eqref{16} we deduce initial condition $u(0)=0$, which is transformed by differential transformation to $U(0)=0$.

First, we apply current approach of using DTM represented by Theorem \ref{t2}.  Transformed equation \eqref{15} has the form
\begin{equation} \label{17}
(k+1) U(k+1) = \frac{1}{a} U(k) - \frac{1}{a} \sum\limits_{l=k}^{N} (-1)^{l-k} \binom{l}{k} a^{l-k} U(l) + a \delta (k).
\end{equation}

If we try to solve \eqref{17} for $N=1$, we get
\begin{equation} \label{18}
U(1) = U(1) +a
\end{equation}
for $k=0$. It implies that there would be a solution only for $a=0$ which is in contradiction to our assumption that $a>0$.

Next, if we solve \eqref{17} for $N=2$, we have
\begin{equation} \label{19}
U(1) = \frac{1}{a} U(1) - \frac{1}{a} \bigl[ U(0) - a U(1) + a^2 U(2) \bigr] + a \delta (0),
\end{equation}
for $k=0$, which after some rearrangements gives
\begin{equation} \label{20}
0=-a U(2) +a,
\end{equation}
hence $U(2)=1$. Calculation for $k=1$ does not provide any new information, while for $k=2$ we get $U(3)=0$. It can be easily computed that for $N \ge 2$, we always get $U(2)=1$ and $U(k) =0$ for $k \ge 3$.

Note that we do not have any information about $U(1)$ for $N \ge 2$. We can interpret it such that $U(1)$ may be an arbitrary constant $C \in \mathbb R$. Taking this into account, we deduce that in this case solution $u(t)$ is in the form
\begin{equation} \label{21}
u(t) = t^2 + C t, \ C \in \mathbb R.
\end{equation}
It is not difficult to verify that all such functions indeed are solutions of \eqref{15} satisfying initial condition $u(0)=0$. It means that using current DTM approach we obtained one-parameter family of solutions of equation \eqref{15} with initial condition $u(0)=0$,  not a unique solution of Cauchy problem \eqref{15}, \eqref{16}, since we did not utilize initial function at all. Furthermore, from \eqref{18} we can see that the result of applying current DTM approach also may depend on the choice of $N$.

On the other hand, applying the new approach we derive relation
\begin{equation} \label{22}
(k+1) U(k+1) = \frac{1}{a} U(k)  - \frac{1}{a} \bigl[ \delta (k-2) -2a \delta (k-1) + a^2 \delta (k) \bigr] + a\delta (k),
\end{equation}
which can be simplified to
\begin{equation} \label{23}
(k+1) U(k+1) = \frac{1}{a} U(k)  - \frac{1}{a} \delta (k-2) +2 \delta (k-1).
\end{equation}
Transformed initial condition is the same as in the previous case, $U(0)=0$. For $k=0,1,2,\dots$ we have
\begin{align}
U(1) &= \frac{1}{a} U(0) &\rightarrow& U(1)=0,\notag \\
2U(2) &= \frac{1}{a} U(1) +2 &\rightarrow& U(2)=1,\notag \\
3U(3) &= \frac{1}{a} U(2) -  \frac{1}{a} &\rightarrow& U(3)=0,\notag \\
 &\vdots& \vdots& \notag
\end{align}
which gives unique solution
\begin{equation} \label{24}
u(t) = t^2, \ t \in [0,a].
\end{equation}
Applying the classical method of steps, we can see that \eqref{24} is unique solution corresponding to Cauchy problem \eqref{15}, \eqref{16} on $[0,a]$, thus the new approach of using DTM is in perfect agreement with well-established results.\\
{\bf Example 3}. Consider delayed differential equation of the third order
\begin{equation}\label{25}
u'''(t) = -u(t) -u(t-0.3) + e ^{-t+0.3}
\end{equation}
subject to the initial function
\begin{equation} \label{26}
\phi (t) =e^{-t}, \ t \leq 0
\end{equation}
and conditions
\begin{eqnarray}\label{27}
u(0)&=& 1, \nonumber \\ 
u'(0)&=& -1, \\
u''(0)&=& 1. \nonumber
\end{eqnarray}
This problem was solved using the Adomian decomposition method (ADM) by Evans and Raslan \cite{evans} and later using current DTM approach by Karakoc and Bereketoglu \cite{karakoc} and again using ADM by Blanco-Cocom et al. \cite{cocom}.

Straightforward observation gives the information that, as Blanco-Cocom et al. \cite{cocom} point out, it is enough to consider only \eqref{25} and \eqref{26}, since conditions \eqref{27} are not independent of initial function $\phi (t)$ defined in \eqref{26}. However, in fact, in all mentioned papers authors did not use initial function \eqref{26} at all, they solved problem \eqref{25}, \eqref{27} which is not a Cauchy problem. Karakoc and Bereketoglu \cite{karakoc} tried to rectify the situation of not using \eqref{26} by excluding this condition from formulation of the studied problem. Unfortunately, this step led to the same curiosity observed in Example 1, when for instance $u(-0.3)$ is not defined. In any of the cases, uniqueness of solution is not guaranteed.

In both papers using ADM the authors obtained approximate solution using iterative scheme containing a triple integral and compared the result to function $u(t) = e^{-t}$ which is a solution of \eqref{25}, \eqref{27} and satisfies \eqref{26} as well.
Karakoc and Bereketoglu \cite{karakoc} solved equation \eqref{25} using current DTM approach without the dependence on the initial function and determined recurrence relation
\begin{equation}\label{28}
(k+1)(k+2)(k+3)U(k+3) = -U(k) -\sum_{h_1=k}^N (-1)^{h_1-k} \binom{h_1}{k} (0.3)^{h_1-k} U(h_1) + \frac{1}{k!} (-1)^k e^{0.3},
\end{equation}
The authors solved \eqref{28} for $N=6,8,10$ and compared obtained approximate solutions to solution $u(t) = e^{-t}$.

In contrast to complicated formulas mentioned above, our approach gives simple recurrence relation
\begin{equation}\label{29}
U(k+3) = \frac{-U(k)}{(k+1)(k+2)(k+3)}, \quad k \geq 0.
\end{equation}
From initial conditions \eqref{27} and recurrence relation \eqref{29}  we have
$$
 U(0)= 1,\ U(1) = -1,\ U(2) = \frac{1}{2}, \ U(3) = \frac{-1}{3!}, \ U(4) = \frac{1}{4!}, \dots , U(k) = \frac{(-1)^k}{k!}, ...
$$
Using inverse differential transformation \eqref{3} we obtain a solution of \eqref{25}, \eqref{26} in the form
\begin{equation} \label{30}
u(t) = \sum_{k=0}^{\infty} \frac{(-1)^k}{k!} t^k = e^{-t}. 
\end{equation}
It is the closed form unique solution of Cauchy problem \eqref{25}, \eqref{26} which cannot be reached using either ADM or current approach of using DTM as described in above mentioned papers \cite{evans}, \cite{cocom} and \cite{karakoc}, only approximation of the solution is possible.

Now consider initial function
\begin{equation}\label{31}
\phi (t) = \frac{1}{2} t^2 -t+1, \ t \in [-0.3,0].
\end{equation}
This function satisfies conditions \eqref{27} as well. However, applying classical method of steps shows that function
\begin{align}
u(t) =& \left( 2.23 + \frac{\e^{0.3}}{3} \right) \e^{-t} + \left( 0.16 - \frac{\e^{0.3}}{3} \right) \e^{\frac{t}{2}} \cos \left( \frac{\sqrt{3}}{2} t \right) + \left( \frac{-0.3}{\sqrt{3}} + \frac{\e^{0.3}}{3 \sqrt{3}} \right) \e^{\frac{t}{2}} \sin \left( \frac{\sqrt{3}}{2} t \right)\notag\\
 &+ \frac{\e^{0.3}}{3} t \e^{-t} - \frac{1}{2} t^2 + 1.3 t - 1.39 \label{32}
\end{align}
is unique solution of Cauchy problem \eqref{25}, \eqref{31} on interval $[0,0.3]$, nevertheless it is completely different from \eqref{30}. Using of our approach leads to recurrence relation
\begin{equation}\label{33}
(k+1)(k+2)(k+3)U(k+3) = -U(k) - \frac{1}{2} \delta (k-2) + 1.3 \ \delta (k-1) - 1.39 \ \delta (k) + \e^{0.3} \frac{(-1)^k}{k!},
\end{equation}
which together with conditions \eqref{27} transformed to $U(0) = 1$, $U(1) = -1$, $U(2) = \frac{1}{2}$ enables to calculate
\begin{alignat}{3}
U(3) &= \frac{-2.39 + \e^{0.3}}{3!}, \ &U(4) &= \frac{2.3 - \e^{0.3}}{4!}, \ &U(5) = \frac{-1 + \e^{0.3}}{5!},\notag\\
U(6) &= \frac{2.39 - 2 \e^{0.3}}{6!}, \ &U(7) &= \frac{-2.3 + 2 \e^{0.3}}{7!}, \ &U(8) = \frac{1 - 2 \e^{0.3}}{8!}, \label{34}\\
U(9) &= \frac{-2.39 + 3 \e^{0.3}}{9!}, &\dots . \notag
\end{alignat}
Precise investigation of this sequence reveals that $U(k)$, $k=0,1,2,\dots$ are coefficients of Taylor expansion of function \eqref{32} in neighbourhood of $0$, thus the calculation resulted in an approximation of unique solution of Cauchy problem \eqref{25}, \eqref{31} on $[0,0.3]$ as expected. Note that in this case, not even approximation of solution is possible using either ADM or current DTM approach as described in papers \cite{evans}, \cite{cocom} and \cite{karakoc}.
\\[5mm]
{\bf Example 4.}
The following Cauchy problem for first order neutral differential equation
\begin{align}
u'(t) + \frac{1}{4} u'(t-1) &= u(t) + u(t-1), \ t \ge 0, \label{r13}\\
u(t) &= -t, \ \text{ for } t \in [-1,0]. \label{r14}
\end{align}
was investigated by Fabiano \cite{fabiano}. The author used semidiscrete approximation scheme to approximate unique solution of problem \eqref{r13}, \eqref{r14}.
The exact solution which can be calculated by method of steps is given by
\begin{equation} \label{r15}
u(t) = t -\frac{1}{4} + \frac{1}{4} \e^t, \ t \in [0,1].
\end{equation}
We remark that the sewing condition
$$
\phi' (0) + \frac{1}{4} \phi' (-1) = \phi (0) + \phi (-1)
$$
 is not fulfilled hence the derivative of solution of problem \eqref{r13}, \eqref{r14} is not continuous at $0$.

Applying DTM combined with method of steps, we obtain relation
\begin{equation} \label{r16}
(k+1) U(k+1) - \frac{1}{4} \delta (k) = U(k) - \delta (k-1) + \delta (k),
\end{equation}
which together with initial condition $u(0)=0$ acquired from \eqref{r14} and transformed to $U(0)=0$ implies
\begin{align} \label{r16}
U(1) - \frac{1}{4} &=  U(0) +1 &\rightarrow& U(1)= \frac{5}{4},\notag \\
2U(2) &= U(1) -1 &\rightarrow& U(2)=\frac{1}{8} = \frac{1}{4 \cdot 2!},\notag \\
3U(3) &= U(2) &\rightarrow& U(3)=\frac{1}{24}=\frac{1}{4 \cdot 3!},\notag \\
4U(4) &= U(3) &\rightarrow& U(4)= \frac{1}{4 \cdot 4!}, \notag \\
 &\vdots& \vdots& \notag
\end{align}
If we write $U(0)$ and $U(1)$ as $U(0)=\frac{1}{4 \cdot 0!} - \frac{1}{4}$ and $U1)=\frac{1}{4 \cdot 1!} +1$ and perform inverse transformation, the solution can be expressed as
\begin{equation} \label{r17}
u(t) = \frac{1}{4 \cdot 0!} - \frac{1}{4} + \left( \frac{1}{4 \cdot 1!} +1 \right) t + \frac{1}{4 \cdot 2!} t^2 + \frac{1}{4 \cdot 3!} t^3 + \frac{1}{4 \cdot 4!} t^4 + \dots + \frac{1}{4 \cdot k!} t^k + \dots = t - \frac{1}{4} + \frac{1}{4} \sum\limits_{k=0}^{\infty} \frac{1}{k!} t^k,
\end{equation}
which is Taylor expansion of exact solution \eqref{r15}. Hence using our approach we are able to identify unique solution of Cauchy problem \eqref{r13}, \eqref{r14} in closed form, which is not the case in Fabiano's paper \cite{fabiano}.

\section*{Conclusion}
\begin{itemize}
\item
We conclude that combination of method of steps and differential transformation method (DTM) presented in this paper is powerful and efficient semi-analytical technique suitable for numerical approximation of a solution of Cauchy problem for wide class of functional differential equations, in particular delayed and neutral differential equations. No discretization, linearization or perturbation is required.
\item
There is no need for calculating multiple integrals or derivatives and less computational work is demanded compared to other popular methods (Adomian decomposition method, variational iteration method, homotopy perturbation method, homotopy analysis method).
\item
Using presented approach, we are able not only to obtain approximate solution, but even there is a possibility to identify unique solution of Cauchy problem in closed form.
\item
A specific advantage of this technique over any purely numerical method is that it offers a smooth, functional form of the solution over a time step.
\item
Another advantage is that using our approach we avoided ambiguities, incorrect formulations and ill-posed problems that occur in recent papers, as we observed in examples.
\item
Finally, a subject of further investigation is to develop the presented technique for equation \eqref{1} with state dependent or time dependent delays.
\end{itemize} 

\section*{Competing interests}
The authors declare that they have no competing interests.

\section*{Authors' contributions}
All authors have made the same contribution. All authors read and approved the final manuscript.

\section*{Acknowledgements}
  \ifthenelse{\boolean{publ}}{\small}{}
The first author was supported by the project CZ.1.07/2.3.00/30.0039 of Brno
University of Technology. The work of the second author was realised in CEITEC - Central European Institute of Technology with research infrastructure supported by the project CZ.1.05/1.1.00/02.0068 financed from European Regional Development Fund and by the project FEKT-S-11-2-921 of Faculty of
Electrical Engineering and Communication, Brno University of Technology.
This support is gratefully acknowledged.

{\ifthenelse{\boolean{publ}}{\footnotesize}{\small}
 \bibliographystyle{bmc_article}  
  \bibliography{bmc_article} }     


\ifthenelse{\boolean{publ}}{\end{multicols}}{}

\end{bmcformat}
\end{document}